\documentclass[12pt,a4paper,reqno]{amsart}
\usepackage{amsfonts}
\usepackage{amsthm}
\usepackage{amsmath}
\usepackage{amscd}
\usepackage[latin2]{inputenc}
\usepackage{t1enc}
\usepackage[mathscr]{eucal}
\usepackage{indentfirst}
\usepackage{graphicx}
\usepackage{graphics}
\usepackage{pict2e}
\usepackage{epic}
\numberwithin{equation}{section}
\usepackage[margin=2.9cm]{geometry}
\usepackage{epstopdf} 

\theoremstyle{plain}
\newtheorem{Th}{Theorem}[section]

\theoremstyle{definition}

\newtheorem{?}[Th]{Problem}

\newcommand{\F}{\mbox{$\mathbb{F}$}}

\begin{document}

\title{Twin prime polynomials in function field}

\author[Sushma Palimar]{Sushma Palimar}

\address{Department of Mathematics\\
Indian Institute of Science\\
Bangalore} 

\email{sushmapalimar@gmail.com, psushma@iisc.ac.in}
 \begin{abstract}
   We resolve  the function field analogue of the conjecture concerning distribution of twin primes in arithmetic progression.

\end{abstract}
\maketitle
\section{Introduction}

It is well known that, if  $\pi(x;q,a)$ denotes the number of primes $p\leq x$ that are $\equiv a\pmod{q}$, and if $(a,q)=1$, then \begin{equation}\pi(x;q,a)\sim\frac{\pi(x)}{\phi(q)}\end{equation} for $x$ suffiently large.  Here $\pi(x)$ is the number of primes $p\leq x$ and $\phi(q)$, the Euler totient function, is the number of residues coprime to $q$.
Motivating question in analytic number theory is the following: for eg., how many prime numbers less than $ x$ satisfy a given condition,  are contained
in a given arithmetic progression $\equiv  a \pmod q?$ The question is trivial when $a$ and $q$ are not relatively prime, but when $(a,q)=1$, one can seek an asymptotic formula for the counting function\begin{equation}
      \pi(x;q,a):= 
      \sum_{\substack{p\leq x\\
      p\equiv a\bmod{q}}}  1.                                                                                                                                                                                          
                                                                                                                                \end{equation}
Computations reveal that the primes up to $x$ are
equi-distributed amongst the arithmetic progressions $mod q,$ once $x$ is little larger than $q$, say $x \geq q^{1+\delta}$
for any fixed $\delta> 0$.
However in many application one wants to know
just how large $x$ needs to be for the primes to be equi-distributed in arithmetic progressions $\pmod{q}.$
The  Bombieri--Vinogradov theorem, is  a major result of analytic number theory, concerning the distribution of primes in arithmetic progressions, averaged over a range of moduli. 
This theorem is a statement about the error term in Dirichlet's theorem  in an arithmetic progression, which
states that primes have a level of distribution 1/2.
 \section{Twin primes in arithmetic progression} The well known conjecture on twin prime states that, there are an infinite number of twin primes and proving this is one of the unsolved problem in number theory. Though substantial progress has been made by Y. Zhang, J. Maynard, T. Tao and others, the problem remains still open.
For a nonzero integer $k$, and for coprime integers $a$, $q$, twin primes in arithmetic progression $a \text{ modulo } q$ is defined  by            
\begin{equation}\label{eq1}
\begin{split}
\psi(x;q,a,2k)=\sum_{\substack{0<m<n\leq x \\ m-n=2k \\ n\equiv a \mathrm{(mod)} q}}\Lambda(m)\Lambda(n)\\
\text{where $\Lambda$ is the von Mangoldt function}
\end{split}
\end{equation}
  
\begin{equation}
\begin{split}
\text{when $(a+2k, q) =1$, 
$\psi$ is asymptotically equal to}\\
H(x;q,2k)={C} \prod_{p\mid qk,p>2}\big(\frac{p-1}{p-2}\big)\frac{x-|2k|}{\phi(q)} \quad \text{ as } x\rightarrow \infty.
\end{split}
\end{equation}
where $$C=(1-\frac{1}{(p-1)^2}\big)$$
 Let $$ E(x;q,a,2k)=\begin{cases}
                  \psi-H \quad if\quad(a+2k,q)=1\\
                  \psi \quad otherwise
                 \end{cases}$$
In most arithmetic applications it is crucial to allow the modulus $q$ to grow  along with $x$. Thus the remainder term in (\ref{eq1}) is of interest.
In literature, studies have been made regarding the distribution of twin primes in arithmetic progressions,
 by A.F. Lavrik \cite{afl}, Maier and  Pomerance \cite{mp} and Hiroshi Mikawa\cite{hm}. Among them, Mikawa's conjectured  asymptote on, ``twin primes in arithmetic progression"  is much closer analog to  Bombieri-Vinogradov asymptote.  A.F.Lavrik \cite{afl} shows, for any $A,B>0$
\begin{equation}\sum_{0<2k\leq x} |E(x;q,a,2k)|\ll x^2(\mathrm{log}x^{-A})\end{equation}
uniformly for $(a,q)=1$ and $q\leq (log x)^{B}$.
Maier and Pomerance \cite{mp} established the inequality
\begin{equation}\sum_{q\leq Q} \max_{(a,q)=1} \sum_{0<2k\leq x}|E(x;q,a,2k)|\ll x^2(\mathrm{log}x)^{-A}\end{equation} and showed that, the above is valid for $Q\leq x^{\delta}$ for some $\delta>0$.
Balog\cite{ab} generalized this to the case of prime multiplets, and extended the range of validity, in the general case, to
$Q\ll x^{1/3}(\log x)^{-B}$ with some $B=B(A)>0$ .
For twin primes in arithmetic progression, Mikawa conjectured the following.
\begin{Th}\label{th1}
Let $A>0 $ be given. There exists $B=B(A)$ such that
\begin{equation}\label{theq} \sum_{q\leq x^{1/2}({log}x)^{-B}} 
\max_{a:(a,q)=1}\sum_{0<2k\leq x}|E(x;q,a,2k)|\ll_{A} x^{2}\mathrm({log}x)^{-A}\end{equation}
where the implied constant depends only on $ A.$
\end{Th}
The asymptote in (\ref{theq}) is obtained by establishing  the following identity, \begin{equation}\label{eq3}
\begin{split}
&\sum_{0<2k\leq x}|E(x;q,a,2k)|^{2}\leq \mathfrak{H}(q)\frac{2}{\phi^{2}(q)}\frac{ x^{3}}{3}\\
\text{where }
&\mathfrak{H}(q)=\prod_{p}(1+\frac{1}{(p-1)^{3}})\prod_{p\mid q}(\frac{(p-1)^{2}}{p^{2}-3p+3}).
\end{split}
\end{equation} 
\section{Twin prime polynomials and some results in function field}
In this short article we resolve the function field versions of conjecture (\ref{eq3}) in arithmetic progression.  The function field version of 
distribution of  prime polynomials in arithmetic progression is done in \cite{KZ} and the following result is derived by Keating and Rudnick. These results are the function field analogoue of conjectures  on distribution of primes by Hooley, Friedlader-Goldston and Montgomery in arithmetic progression. %\cite{hml}.
\begin{Th}
   \begin{enumerate}
    \item Given a finite field $\F_{q}$, let $Q\in \F_{q}[t]$ be a polynomial of  positive degree, and $1\leq n\leq\ \mathrm{deg} Q$. 
    Then
   $ G(n;Q)=nq^{n}-\frac{q^{2n}}{\Phi(Q)}+O(n^2 q^{\frac{n}{2}})+O((deg Q)^2)$
    \item Fix $n\geq 2$. Given a sequence of finite fields $\F_{q}$ and square-free polynomials $Q(T)\in \F_{q}[T]$ of positive degree with 
    $n\geq \mathrm{deg}Q-1,$ then as $q\rightarrow \infty$,\[G(n;Q)\sim 
    q^{n}(deg Q-1).\]
    \end{enumerate}
 \end{Th}

Let $\F_{q}[t]$ be the ring of polynomials over the finite field $\F_{q}$ with $q$ elements, $q=p^{\nu}, p: \text{ prime}.$
 Let  $\mathcal{P}_{n}=\{f\in \F_{q}[t]| \mathrm{deg} f=n \}$ be the set of all polynomials of degree $n$ and 
 $\mathcal{M}_{n}\subset\mathcal{P}_{n}$ be the subset of monic polynomials of degree $n$ over $\F_{q}$.
 The von Mangoldt function in this case is defined as 
  $$\Lambda(N)=\begin{cases}\mathrm{ deg}P \text{ if } N=cP^{k}, P \text{ is an irreducible monic polynomial and }
  c\in \F_{q}^{\times}\\
  0 \text{ otherwise}
              \end{cases}$$
 
 We study the function field analog of distribution of twin prime polynomials in arithmetic progressions.
The polynomial version of twin prime conjecture was proved by Bender and Pollack in \cite{AP}.
Let $A$ in $\F_{q}[t]$ be a polynomial of $\mathrm{deg} K <n$. Then the number of monic irreducible polynomials of degree $n$ for which both $F$ and $F+K$ are irreducible is 
$$\pi_{2}(n;q)\approx\frac{q^{n}}{n^{2}}$$ Here $q$ is large as compared with $n^2$.
It follows from the work of Bary-Soroker  and Bender and Pollack,
that for fixed $n$, in the limit of large (odd) $q$, and given a polynomial $K\neq 0$ and $deg K<n$, the analogue of the Hardy-Littlewood conjecture for twin primes holds in the form
\begin{equation}\label{asym}
 \sum_{f\in \mathcal{M}_n}\Lambda(f)\Lambda(f+K)\sim q^{n}+O_n(q^{n-1/2})
\end{equation}
For $r$ poynomials $h_1,..,h_r$ of degree less than $n$, Hardy-Littlewood prime tuple conjecture for the rational function field is
\begin{equation}\label{sminsert} \sum_{\substack{f\in \F_{q}[t]\\deg f> deg h_i}}\Lambda(f+h_1)...\Lambda(f+h_r)\sim \mathfrak{C}(h_1,...h_r) q^{n} \text{ as }
 q^n\rightarrow \infty.\end{equation}
%%%%%%%%%%%%%%%%%%%%%
where the singular series $\mathfrak{C}(h_1,...h_r)$ is given by
\begin{equation}\mathfrak{C}(h_1,...h_r)=\prod_{P}(1-\frac{1}{p})^{-r}
(1-\frac{\nu_{h_1,h_2,...,h_r}}{|P|})\end{equation}
with the product involving all monic irreducibles $P$ and 
\begin{equation}\label{asym1}
\begin{split}
 &\nu_{h_1,h_2,...,h_r}=\#\{A \bmod{P}:(A+h_1)...(A+h_r)\equiv 0\bmod{P}\}\\
Note \, that,\\
&\mathfrak{C}({h_1,h_2,...,h_r})=1+O(\frac{1}{q})
\end{split}
\end{equation}
Finding irreducibles in prime tuples, which inturn leads to counting of rational points of an absolutely irreducible variety  over a finite field $\F_q$. 
 Then the required asymptotic follows by applying Lang - Weil estimate. 
\subsection{Arithmetic Progression}
For a polynomial $Q\in\F_{q}[t]$ of positive degree, and $A\in\F_{q}[t]$ coprime to $Q$, $\mathrm{deg}K<\mathrm{deg}n$ and any $n>0,$ set
\begin{equation}\label{apeq}\psi(n;Q,A,K)=\sum_{\substack {f\in \mathcal{M}_{n}\\f\equiv A\mathrm{mod}Q}}\Lambda(f)\Lambda(f+K)\end{equation}
We know that, \begin{equation}
               \sum_{f\in \mathcal{M}_{n}}\Lambda(f)\Lambda(f+K) \sim q^{n}
              \end{equation}
Twin prime polynomial theorem in arithmetic progression gives the following asymptote \cite{jpkz, KZ}.
\begin{equation}
\begin{split}
\label{twinp} \psi(n;Q,A,K)= & \sum_{\substack{f\in \mathcal{M}_{n}\\ f\equiv A\mathrm{mod}Q }}\Lambda(f)\Lambda(f+K)\\                                                                                                                                                                                                                                                                             & \approx \frac{q^{n-\mathrm{deg}Q}|Q|}{\Phi(Q)} \approx \frac{q^{n}}{\Phi(Q)}
\end{split}\end{equation}
where $\Phi(Q)$ is the Euler totient function, namely number of reduced residues modulo $Q$.
We wish to show an analog of conjecture \ref{eq3} in the limit of large finite field size, that is $q\rightarrow \infty.$ Define,
\begin{equation}
E(n;Q,A,K)=|\psi(n;Q,A,K)-\frac{q^{n}}{\Phi(Q)}|\end{equation}
\begin{equation}\label{evaleq}
\begin{split}
&H(n;Q):=\sum_{\substack{A \mathrm{mod}  Q\\ gcd(A,Q)=1}}|\psi(n;Q,A,K)-\frac{q^{n}}{\Phi(Q)}|^2\\
&{ and }\\
&G(n,q)=\sum_{\substack{A \mathrm{mod}  Q\\ gcd(A,Q)=1}}|\psi(n;Q,A,K)-\frac{q^{n}}{\Phi(Q)}|
\end{split}
\end{equation}
We have the analogue of function field results.
\begin{Th}\label{an-kz}
 Fix $n\geq  \mathrm{deg}Q-1.$ In the limit $q\rightarrow\infty$,
 \[H(n;Q)=n^2q^{2n} \text{ and } G(n;q)\sim nq^{n}\] \end{Th}
 \begin{Th}\label{an-fun}
For $n\geq \mathrm{deg}Q-1$  and $0\neq K<n,$ we have,
\begin{equation}\label{eqth}
 \sum_{k=0}^{n-1}\sum_{K\in \mathcal{M}_{k}}E(n;Q,A,K)^2\approx  
 \frac{q^{3n}}{\phi(q)^2}\end{equation}
\end{Th}
We can compare our result in (\ref{eqth}) with that of (\ref{eq3})
 when the following correspondence is used.
 \[Q\leftrightarrow ||Q||=q^{\mathrm{deg}Q}, \quad \mathrm{log}Q \leftrightarrow \mathrm{deg }Q, \quad X\leftrightarrow q^{n},\quad  \mathrm{log}X\leftrightarrow n\]
 \[\]
 \subsection{Proof of theorem \ref{an-kz}}
 We start evaluating equation (\ref{evaleq})
\begin{equation}\label{sub-var}
 H(n;Q)=\sum_{gcd(A,Q)=1}\psi(n;Q,A)^{2}-2\frac{q^{n}}{\Phi(Q)}\sum_{gcd(A,Q)=1}\psi(n;Q,A)+\frac{q^2n}{\Phi(Q)}
\end{equation}
The first moment of $\psi(n;Q,A,K)$ is 
\begin{equation}
\sum_{gcd(A,Q)=1}\psi(n;Q,A,k)=\sum_{\substack{\mathrm{deg}f=n\\\mathrm{deg}K<n\\gcd(f,Q)=1}}\Lambda(f)\Lambda(f+K)\end{equation}
\begin{equation}\label{subq1}
=\sum_{\substack{deg f=n\\deg K<n}}\Lambda(f)\Lambda(f+K)-\sum_{\substack{\mathrm{deg}f=n\\\mathrm{deg}K<n\\gcd(f,Q)>0}}\Lambda(f)\Lambda(f+K)
=q^{n}-\text{ terms of lower order }
\end{equation}
%We can compare the result if we make the following dictionary.
For the second moment of $\psi(n;Q,A,K)$, we have
\begin{equation}
\sum_{gcd(A,Q)=1}\psi(n,Q,A,K)^{2}=\sum_{\substack{\mathrm{deg}f=\mathrm{deg}g=n\\f\equiv g\bmod{Q}\\gcd(f,Q)=1}
}\Lambda(f)\Lambda(f+K_1)\Lambda(g)\Lambda(g+K_2)
\end{equation}
\begin{equation}\label{splitap}
=\sum_{\substack{\mathrm{deg}f=\mathrm{deg}g=n\\f\equiv g\bmod{Q}\\ f\neq g\\gcd(f,Q)=1}
}\Lambda(f)\Lambda(f+K_1)\Lambda(g)\Lambda(g+K_2)+  \sum_{\substack{deg f=n\\
gcd(f,Q)=1\\\mathrm{deg}(K)<n}}\Lambda(f)^{2}\Lambda(f+K)^{2}
\end{equation}
Consider the first sum
\begin{equation}\label{val1}
 \sum_{\substack{deg f=deg g=n\\
f\equiv g\bmod Q\\f\neq g\\
gcd(f,Q)=1}}  \Lambda(f)\Lambda(g)\Lambda(f+K_1)\Lambda(g+K_2)    
\end{equation}
Since $f,g\in \mathcal{M}_{n}$ we can write $f\equiv g \bmod Q 
\text{  as } g=f+JQ, \text{ for} J\neq 0 \text{ and }\mathrm{deg}J<n-\mathrm{deg} Q$.
Hence, the above equation is
\begin{equation}\label{eval}       
\begin{split}
\sum_{\substack{deg f=deg g=n\\f\equiv g\bmod Q\\f\neq g\\gcd(f,Q)=1}} 
\Lambda(f)\Lambda(f+JQ)\Lambda(f+K_1)\Lambda(j+JQ+K_2) \text{ as}
  &\sum_{j=0}^{n-deg Q}\sum_{degJ<n-deg Q}                                                                                       \psi_{fg}(n;JQ) \\
  \psi_f(n;K)=\sum_{\substack{\mathrm{deg} f= n\\ deg K<n}}\Lambda(f)\Lambda(f+K)\approx
q^{n}(\mathfrak{C}(K))\\
\text{where } \mathfrak{C}(K) \text{ is given by  (\ref{asym1}) and }
 \mathfrak{C}(K)=1+\frac{1}{q} 
\end{split}
\end{equation}
%\[\sim q^{n}(q-1)\sum_{j=0}^{n-deg Q}\sum_{degJ<n-deg Q}\]
  So, for fixed $n$ as $q\rightarrow \infty,$ from Hardy -Littlewood prime tuple conjecture for rational function fields, \cite{LB}
 \begin{equation}\label{subq}
  \psi_{fg}\sim {q^{n}(1+\frac{1}{q})}
 \end{equation}
To evaluate equation (\ref{val1}) we can restrict the $J$-sum in equation (\ref{eval})
to monics, by multiplying by $(q-1)$.  Applying the estimate (\ref{asym}) and  in (\ref{eval}), we have 
\begin{equation}
\begin{split}\label{mainsub1}
 & \sum_{\substack{deg f=deg g=n\\
f\equiv g\bmod Q\\f\neq g\\
gcd(f,Q)=1}}  \Lambda(f)\Lambda(g)\Lambda(f+K)\Lambda(g+K)\\
&\approx q^{n}(q-1)\sum_{j=0}^{n-\mathrm{deg} Q}\sum_{\substack{\mathrm{deg} J=j\\J \neq 0\\
J \text{ monic }}} \mathfrak{C}(JQ) \text{ as } q^{n}\rightarrow\infty.
\end{split}
\end{equation}
\text{The estimate for } $\mathfrak{C}(JQ)$                                          \text{ is  given in \cite{KZ}, Appendix A}
\begin{equation}\label{Jsum}
\sum_{\substack{\mathrm{deg} J=j\\J \neq 0\\
J \text{ monic }}} \mathfrak{C}(JQ)\sim \frac{q^{j}|Q|}{\phi(Q)}-\frac{1}{1-q}\end{equation}
Hence,  the right hand side of (\ref{mainsub1}) becomes
\begin{equation}\label{subeq1} \sum_{\substack{deg f=n\\deg g=n\\ f\equiv g \bmod Q\\ f\neq g \\deg K<n}}\Lambda(f)\Lambda(g)\Lambda(f+K)\Lambda(g+K)\sim \frac{q^{2n}}{\phi(Q)}\end{equation}
When $f=g$, we have 
\begin{equation}\label{lamsq}
\sum_{\substack{deg f=n\\
gcd(f,Q)=1\\\mathrm{deg}(K)<n}}\Lambda(f)^{2}\Lambda(f+K)^{2}\sim n^{2}q^{2n}\end{equation}
This follows from, [\cite{jpkz}, Lemma 3.1],
\begin{equation}\label{sqsum}
\sum_{f\in \mathcal{M}_{n}}\Lambda(f)^{2}=nq^{n}+O(n^{2}q^{n/2})\end{equation}
Therefore, substituting, (\ref{subq1}),  (\ref{subq}) and (\ref{subeq1}), in   (\ref{sub-var}), we have
\begin{equation}\label{fineq}
\begin{split}
 &H(n;Q)=\sum_{gcd(A,Q)=1}\psi(n;Q,A)^{2}-2\frac{q^{n}}{\Phi(Q)}\sum_{gcd(A,Q)=1}\psi(n;Q,A)+\frac{q^2n}{\Phi(Q)}\\
  &H(n;Q)\approx n^{2}q^{2n} \text{ and } \\
  &G(n;Q)=\sum_{gcd(A,Q)=1} |\psi(n;Q,A)-\frac{q^{n}}{\Phi(Q)}|\approx nq^{n}\\
  \end{split}
 \end{equation}
 which proves  estimate in theorem \ref{an-kz}.
 To prove Theorem, \ref{an-fun},    we note that,
 \[E(n;Q,A,K)=|\psi(n;Q,A,K)-\frac{q^{n}}{\Phi(Q)}|\]
  Evaluation equation ({\ref{eqth}}), Since,
  \[E(n;Q,A,K)^{2}=\psi(n;Q,A)^{2}-2\frac{q^{n}}{\Phi(Q)} \psi(n;Q,A)+\frac{q^{2n}}{\Phi(Q)^2}\]
Since, \[E(n;Q,A,K)^{2}=\psi(n;Q,A)^{2}-2\frac{q^{n}}{\Phi(Q)} \psi(n;Q,A)+\frac{q^{2n}}{\Phi(Q)^2}\]
and
\begin{equation}\label{subin}
\sum_{k=0}^{n-1}\sum_{K\in \mathcal{M}_{k}}E(n;Q,A,K)^2=
  \sum_{k=0}^{n-1}\sum_{K\in \mathcal{M}_{k}}\big(\psi(n;Q,A)^{2}-2\frac{q^{n}}{\Phi(Q)} \psi(n;Q,A)+\frac{q^{2n}}{\Phi(Q)^2}\big)\end{equation}
 
\begin{equation}\label{ersub1} \sum_{k=0}^{n-1}\sum_{K\in \mathcal{M}_{k}}\psi(n;Q,A)=\sum_{k=0}^{n-1}\sum_{K\in \mathcal{M}_k}\sum_{\substack{\mathrm{deg}f=n\\\mathrm{deg}K<n\\gcd(f,Q)=1}}\Lambda(f)\Lambda(f+K)\end{equation}
 R.H.S. of (\ref{ersub1}) is
 \begin{equation}\label{rhsersub1}
 \frac{q^{n}}{\phi(Q)}\sum_{k=0}^{n-1}\sum_{K\in \mathcal{M}_k} \mathfrak{C}(K)
 \quad \text{ and }\quad \sum_{K\in \mathcal{M}_k} \mathfrak{C}(K)\sim q^{k}\end{equation}
Substituting (\ref{rhsersub1}) in (\ref{ersub1}), we have
\begin{equation}\label{sub2}
 \sum_{k=0}^{n-1}\sum_{K\in \mathcal{M}_{k}}\psi(n;Q,A)\sim\frac{q^{2n}}{\phi(Q)}
\end{equation}
To calculate the second moment we have from (\ref{splitap})
\[\sum_{k=0}^{n-1}\sum_{K\in \mathcal{M}_{k}}\big(\psi(n;Q,A)^{2}=\]
\begin{equation}\label{ap1}
\sum_{k=0}^{n-1}\sum_{K\in \mathcal{M}_{k}}\sum_{\substack{\mathrm{deg}f=n\\\mathrm{deg}K<n\\f\equiv A\bmod{q}}}\Lambda(f)\Lambda(f+K_1)\Lambda(g)\Lambda(g+K_2)+
\sum_{k=0}^{n-1}\sum_{K\in \mathcal{M}_{k}}\sum_{\substack{\mathrm{deg}f=n\\\mathrm{deg}K<n\\f\equiv A\bmod{q}}}\Lambda(f)^{2}\Lambda(f+K)^{2}
\end{equation}
from (\cite{KZ}, A.11), when $f=g$, 
\[\sum_{(A,Q)=1}\psi(n;Q,A)^{2}=\sum_{\substack{\mathrm{deg}f=n\\gcd(f,Q)=1}}\Lambda(f)^{2}\approx nq^{n}+O(n^2q^{n/2})\]
And, we have from (\ref{lamsq})
\begin{equation}
 \sum_{\substack{\mathrm{deg}f=n\\\mathrm{deg}K<n\\gcd(f,Q)=1}}\Lambda(f)^{2}\Lambda(f+K)^{2}\approx n^2q^{2n}
\end{equation}
Clearly,
\begin{equation}\label{sub1}
  \sum_{k=0}^{n-1}\sum_{K\in \mathcal{M}_{k}}\sum_{\substack{\mathrm{deg}f=n\\\mathrm{deg}K<n\\f\equiv A\bmod{q}}}\Lambda(f)^{2}\Lambda(f+K)^{2}\approx \frac{1}{\Phi(Q)^{2}}n^2q^{3n}
\end{equation}
from (\ref{subeq1}), we have
\[\sum_{\substack{deg f=n\\deg g=n\\ f\equiv g \bmod Q\\(f,Q)=1\\ f\neq g \\deg K<n}}\Lambda(f)\Lambda(g)\Lambda(f+K)\Lambda(g+K)\sim \frac{q^{2n}}{\Phi(Q)^2}\]
So now we arrive at,
\begin{equation}\label{sub11}
\sum_{k=0}^{n-1}\sum_{K\in \mathcal{M}_{k}}\sum_{\substack{deg f=n\\deg g=n\\ f\equiv g \bmod Q\\ f\neq g \\deg K<n}}\Lambda(f)\Lambda(g)\Lambda(f+K)\Lambda(g+K)\approx \frac{q^{3n}}{\Phi(Q)^2}\end{equation}

Substituting (\ref{sub1}), (\ref{sub11}) and (\ref{sub2}) in(\ref{subin}), we obtain the desired 
\[\sum_{k=0}^{n-1}\sum_{K\in \mathcal{M}_{k}}E(n;Q,A,K)^2\approx \frac{q^{3n}}{\Phi(Q)^2}\]
as required, and proof of Theorem (\ref{an-fun}) is complete. Also, same techniques can be used to find, the  distribution of twin primes in short intervals.
\section{Acknowledgement}
 The author would like to express her sincere gratitude to  Prof. Soumya Das, Department of Mathematics, IISc, Bangalore, for suggesting the problems, and  for lot more helpful discussions  while preparing this manuscript.
  This research work was supported by the Department of Science and Technology, Government of India,  under the
 Women Scientist Scheme [SR/WOS-A/PM-31-2016].
 
\bibliographystyle{plain}

\end{document}